\documentclass[11pt,a4paper]{article}

\usepackage{epic,eepic,epsfig,amsmath,amssymb,amsthm}
\usepackage{t1enc}

\usepackage[francais,english]{babel}

\usepackage{caption}

\usepackage{color}

\usepackage{bbm}

\def\virgp{\raise 2pt\hbox{,}}

\renewcommand{\geq}{\geqslant}
\renewcommand{\leq}{\leqslant}
\def\N{{\mathbb N}}

\def\R{{\mathbb R}}

\def\virgp{\raise 2pt\hbox{,}}
\def\cdotpv{\raise 2pt\hbox{;}}

\def\1{\mathbbm{1}}
\newcommand{\ds}{\displaystyle}

\newtheorem{theorem}{Theorem}[section]

\newtheorem{proposition}[theorem]{Proposition}

\newtheorem{pte}[theorem]{Property}

\theoremstyle{remark}
\newtheorem{remark}{Remark}[section]

\theoremstyle{definition}
\newtheorem{definition}{Definition}[section]

\newtheorem*{notation}{Notation}

\theoremstyle{definition}

\theoremstyle{definition}

\begin{document}

\title{Laplacian, on the Sierpi\'{n}ski tetrahedron}

\author{Nizare Riane, Claire David}

\maketitle
\centerline{Sorbonne Universit\'es, UPMC Univ Paris 06}

\centerline{CNRS, UMR 7598, Laboratoire Jacques-Louis Lions, 4, place Jussieu 75005, Paris, France}

\begin{abstract}
Numerous work revolve around the Sierpi\'nski gasket. Its three-dimensional analogue, the Sierpi\'nski tetrahedron~$\mathfrak {ST}$. In~\cite{FengSong}, the authors discuss the existence of the Laplacian on~$\mathfrak {ST}$. In this work, We go further and give we give the explicit spectrum of the Laplacian, with a detailed study of the first eigenvalues. This enables us to obtain an estimate of the spectral counting function, by applying the results given in~\cite{Kigami1998},~\cite{Strichartz2012}.
\end{abstract}

\maketitle
\vskip 1cm

\noindent \textbf{Keywords}: Laplacian - Sierpi\'nski tetrahedron.

\vskip 1cm

\noindent \textbf{AMS Classification}:  37F20- 28A80-05C63.
\vskip 1cm

\vskip 1cm

\section{Introduction}

The Laplacian plays a major role in the mathematical analysis of partial differential equations. Recently, the work of J. Kigami ~\cite{Kigami1989}, ~\cite{Kigami1993}, taken up by
 R.~S.~Strichartz~\cite{Strichartz1999},~\cite{KigamiStrichartzWalker2001}, allowed the construction of an operator of the same nature, defined locally, on graphs having a fractal character: the triangle of Sierpi\'nski, the carpet of Sierpi\'nski, the diamond fractal, the Julia sets, the fern of Barnsley. \\

J.~Kigami starts from the definition of the Laplacian on the unit segment of the real line. For a  double-differentiable function~$u$ on~$ [0,1] $, the Laplacian~ $ \Delta \, u $ is obtained as a second derivative of~$ u $ on~$ [0,1] $. For any pair~$ (u, v) $ belonging to the space of functions that are differentiable on $ [0,1] $, such that:
$$v(0)=v(1)=0$$

\noindent he puts the light on the fact that, taking into account:

$$  \displaystyle \int_0^1 \left (\Delta u \right)(x)\,v(x)\,dx=-   \displaystyle \int_0^1 u'(x)\,v'(x)\,dx
=- \displaystyle \lim_{n \to + \infty} \sum_{k=1}^n  \displaystyle \int_{ \frac{k-1}{n}}^{ \frac{k }{n}} u'(x)\,v'(x)\,dx $$

\noindent if~$\varepsilon > 0$, the continuity of~$u'$ and~$v'$ shows the existence of a natural rank~$n_0$ such that, for any integer~$n \geq n_0$, and any real number~$x$ of~$\left [\displaystyle \frac{k-1}{n}, \displaystyle\frac{k }{n} \right]$,~$1 \leq k \leq n$:

$$\left | u'  (x)- \displaystyle   \frac{ u \left (\displaystyle \frac{k}{n} \right) -u \left (\displaystyle \frac{k-1}{n} \right)}{\frac{1}{n} }\right| \leq \varepsilon
\quad , \quad \left | v'  (x)- \displaystyle   \frac{ v \left (\displaystyle \frac{k}{n} \right) -v \left (\displaystyle \frac{k-1}{n} \right)}{\frac{1}{n} } \right| \leq \varepsilon$$
\noindent the relation:

$$  \displaystyle \int_0^1 \left (\Delta u \right)(x)\,v(x)\,dx=
-\displaystyle \lim_{n \to + \infty} n \,\displaystyle \sum_{k=1}^n \left (u \left (\displaystyle \frac{k}{n} \right) -u \left (\displaystyle \frac{k-1}{n} \right) \right) \,
\left (v \left (\displaystyle \frac{k}{n} \right) -v \left (\displaystyle \frac{k-1}{n} \right) \right)
 $$

\noindent enables one to define, under a weak form, the Laplacian of~$u$, while avoiding first derivatives. It thus opens the door to Laplacians on fractal domains.\\

Concretely, the weak formulation is obtained by means of Dirichlet forms, built by induction on a sequence of graphs that converges towards the considered domain. For a continuous function on this domain, its Laplacian is obtained as the renormalized limit of the sequence of graph Laplacians.  \\

Numerous work revolve around the Sierpi\'nski gasket. Its three-dimensional analogue, the Sierpi\'nski tetrahedron~$\mathfrak {ST}$, obtained by means of an iterative process which consists in repeatedly contracting a regular~$3-$simplex to one half of its original height, put together four copies, the frontier corners of which coincide with the initial simplex, appears as a natural extension. Yet, very few works concern~$\mathfrak {ST}$ in the existing literature. In~\cite{FengSong}, the authors discuss the existence of the Laplacian on~$\mathfrak {ST}$. Yet, they do not give what appears to be of the higher importance, i.e. the spectrum of the Laplacian. In~\cite{Hernandez}, generalizations of the Sierpi\'nski gasket to higher dimensions are considered. Yet, despite interesting results, there are a few mistakes, and no study at all of the Dirichlet forms, whereas they are the  obligatory passage to the determination of the Laplacian. \\

We go further and, after a detailed study, we give the explicit spectrum of the Laplacian, with a specific presentation of the first eigenvalues. This enables us to obtain an estimate of the spectral counting function (analogous of Weyl's law), by applying the results given in~\cite{Kigami1998},~\cite{Strichartz2012}.

 \begin{figure}[h!]
 \center{\psfig{height=8cm,width=10cm,angle=0,file=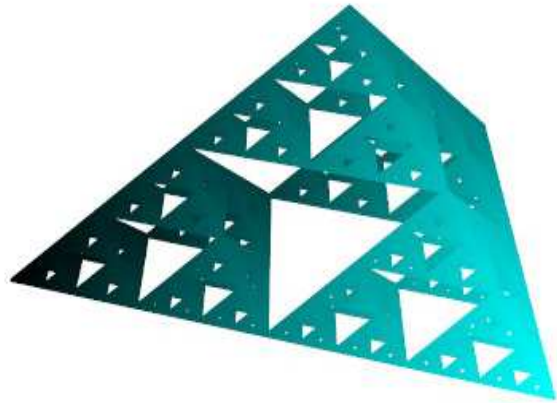}}\\
\caption{The Sierpi\'{n}ski tetrahedron.}
\label{fig1}
\end{figure}

 \begin{figure}[h!]
 \center{\psfig{height=8cm,width=10cm,angle=0,file=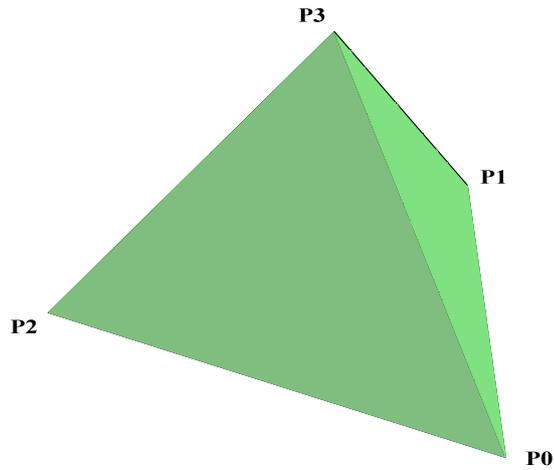}}\\
\caption{The initial tetrahedron.}
\end{figure}

 \begin{figure}[h!]
 \center{\psfig{height=8cm,width=10cm,angle=0,file=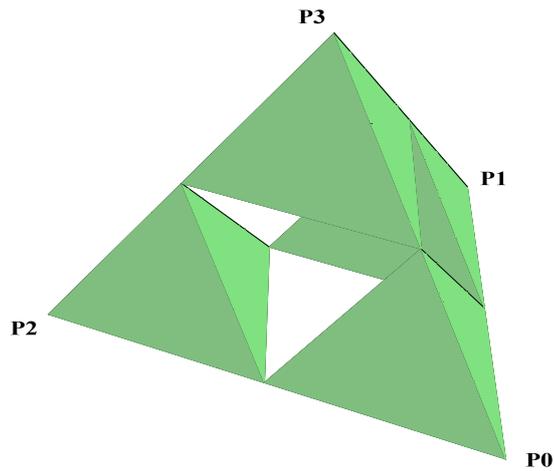}}\\
\caption{The tetrahedron after one iteration.}
\end{figure}

 \begin{figure}[h!]
 \center{\psfig{height=8cm,width=10cm,angle=0,file=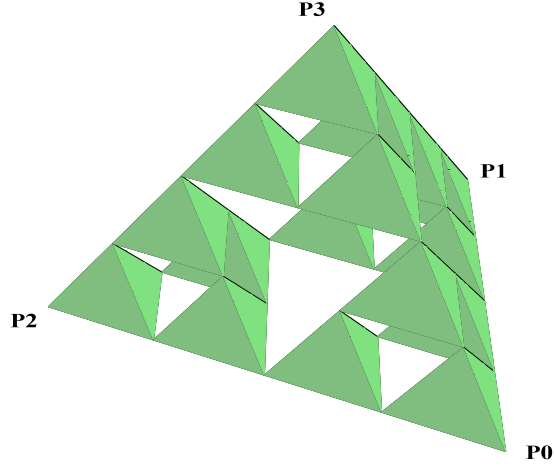}}\\
\caption{The tetrahedron after two iterations.}
\end{figure}

The Sierpi\'{n}ski tetrahedron is a self-similar set which has many beautiful properties.

\section{Self-similar Sierpi\'{n}ski tetrahedron}

We place ourselves, in the following, in the euclidian space of dimension~3, referred to a direct orthonormal frame. The usual Cartesian coordinates are~$(x,y,z)$.\\
\noindent  Let us denote by~$P_0$,~$P_1$,~$P_2$,~$P_3$ the points:

$$P_0= (0,0,0) \quad , \quad P_1= (6^{\frac{1}{3}},0,6^{\frac{1}{3}}) \quad , \quad P_2= (0,6^{\frac{1}{3}},6^{\frac{1}{3}}) \quad , \quad P_3= (0,0,6^{\frac{1}{3}}) $$

\vskip 1cm

Let us introduce the iterated function system of the family of~maps from~$\R^2$ to~$\R^2 $:
$$\left \lbrace f_{0},...,f_{3} \right \rbrace$$
\noindent where, for any integer~$i$ belonging to~\mbox{$\left \lbrace 0,...,3  \right \rbrace$}, and any~$X \,\in\R^2$:
$$f_i(X)=\displaystyle \frac{X+P_i}{2}$$

\begin{remark}
\noindent The family~$\left \lbrace f_{0},...,f_{3} \right \rbrace$ is a family of contractions from~$\R^2$ to~$\R^2$, the ratio of which is:
$$\displaystyle\frac{1}{2}$$

 \noindent For any integer~$i$ belonging to~\mbox{$\left \lbrace 0,\hdots,3  \right \rbrace$},~$P_i$ is the fixed point of~$f_i$.

\end{remark}

\vskip 1cm

\begin{pte}

According to~\cite{Hutchinson1981}, there exists a unique subset $\mathfrak{ST} \subset \R^3$ such that:
\[\mathfrak{ST} = \underset{  i=0}{\overset{3}{\bigcup}}\, f_i(\mathfrak{ST})\]
\noindent which will be called the Sierpi\'{n}ski tetrahedron.\\

\noindent For the sake of simplicity, we set:

$$F=\underset{  i=0}{\overset{3}{\bigcup}}\, f_i  $$

\end{pte}

\vskip 1cm

\begin{definition}\textbf{Hausdorff dimension of the Sierpi\'{n}ski tetrahedron~$\mathfrak{ST}$}\\
 \noindent The Hausdorff dimension of the Sierpi\'{n}ski tetrahedron~$\mathfrak{ST}$ is:

$$D_{\mathfrak {ST}}=\ln_2 4=\displaystyle \frac{\ln 4 }{\ln 2 }=2$$
\end{definition}

\vskip 1cm

\begin{definition}

\noindent We will denote by~$V_0$ the ordered set, of the points:

$$\left \lbrace P_{0},...,P_{3}\right \rbrace$$

\noindent The set of points~$V_0$, where, for any~$i$ of~\mbox{$\left \lbrace  0,1,2  \right \rbrace$}, the point~$P_i$ is linked to the point~$P_{i+1}$, constitutes an oriented graph, that we will denote by~$ {\mathfrak {ST}}_0$.~$V_0$ is called the set of vertices of the graph~$ {\mathfrak {ST}}_0$.\\

\noindent For any natural integer~$m$, we set:
$$V_m =\underset{  i=0}{\overset{3}{\bigcup}}\, f_i \left (V_{m-1}\right )$$

\noindent The set of points~$V_m$, where two consecutive points are linked, is an oriented graph, which we will denote by~$ {\mathfrak {ST}}_m$.~$V_m$ is called the set of vertices of the graph~$ {\mathfrak {ST}}_m$. We will denote, in the following, by~${\cal N}_m$ the number of vertices of the graph~$ {\mathfrak {ST}}_m$.

\end{definition}

\vskip 1cm
\begin{pte}

The set of vertices~$\left (V_m \right)_{m \in\N}$ is dense in~$ {\mathfrak {ST}}$.

\end{pte}

\vskip 1cm
\begin{proposition}
Given a natural integer~$m$, we will denote by~$\mathcal{N}_m$ the number of vertices of the graph ${\mathfrak {ST}}_m$. One has then, for any pair of integers~\mbox{$(i,j)\,\in\,\left \lbrace 0, \hdots, 3 \right \rbrace ^2$}:
\[f_i(P_j)=f_j(P_i)\]
and
\[\mathcal{N}_m=4\times \mathcal{N}_{m-1}-6\]\\
\end{proposition}

\vskip 1cm
\begin{proof}
Let us note that, for any pair of integers~$(i,j) \,\in \,\{0,1,2,3\}^2$:
\[f_i(P_j)=\frac{1}{2}\left(P_j+P_i \right)=f_j(P_i)\]
\noindent Thus, six points appear twice in ${\mathfrak {ST}}_m$.\\
\noindent The second point results from the fact that  the graph~${\mathfrak {ST}}_m$ is obtained by applying four similarities to~${\mathfrak {ST}}_{m-1}$.\\
\end{proof}

 \vskip 1cm

\begin{definition}\textbf{Consecutive vertices of~${\mathfrak {ST}}$ }\\

\noindent Two points~$X$ and~$Y$ of~${\mathfrak {ST}}$ will be called \textbf{\emph{consecutive vertices}} of~${\mathfrak {ST}} $ if there exists a natural integer~$m$, and an integer~$j $ of~\mbox{$\left \lbrace  0,1,2,3 \right \rbrace$}, such that:

$$X = \left (f_{i_1}\circ \hdots \circ f_{i_m}\right)(P_j) \quad \text{and} \quad Y = \left (f_{i_1}\circ \hdots \circ f_{i_m}\right)(P_{j+1})
\qquad\left \lbrace i_1,\hdots, i_m \right \rbrace \,\in\,\left \lbrace  0,1,2,3 \right \rbrace^m $$

\end{definition}

\vskip 1cm

\begin{definition}\textbf{Polyhedral domain delimited by~$ {\mathfrak {ST}}_m  $,~$m\,\in\,\N $}\\

\noindent For any natural integer~$m$, well call \textbf{polyhedral domain delimited by~$ {\mathfrak {ST}}_m  $},
and denote by~\mbox{$ {\cal D} \left ( {{\mathfrak {ST}}_m}\right) $}, the reunion of the~$4^m$ tetrahedra of~${\mathfrak {ST}}_m$.
\end{definition}
\vskip 1cm

\begin{definition}\textbf{Polyhedral domain delimited by the Sierpi\'{n}ski Tetrahedron~$  {\mathfrak {ST}}$ }\\

\noindent We will call \textbf{polyhedral domain delimited by~$  {\mathfrak {ST}}$}, and denote by~\mbox{$ {\cal D} \left ( {\mathfrak {ST}}\right) $}, the limit:
$$ {\cal D} \left ( {\mathfrak {ST}}\right)  = \displaystyle \lim_{n \to + \infty} {\cal D} \left ({\mathfrak {ST}}_m \right) $$

\end{definition}

\begin{definition}\textbf{Word, on ~${\mathfrak {ST}}$}\\

\noindent Let~$m  $ be a strictly positive integer. We will call \textbf{number-letter} any integer~${\cal W}_i$ of~\mbox{$\left \lbrace 0, 1,2,3 \right \rbrace $}, and \textbf{word of length~$|{\cal W}|=m$}, on the graph~$\mathfrak {ST}$, any set of number-letters of the form:

$${\cal W}=\left ( {\cal W}_1, \hdots, {\cal W}_m\right)$$

\noindent We will write:

$$f_{\cal W}= f_{{\cal W}_1} \circ \hdots \circ  f_{{\cal W}_m}  $$

\end{definition}

 \vskip 1cm
\begin{definition}\textbf{Edge relation, on~${\mathfrak {ST}}$}\\

\noindent Given a natural integer~$m$, two points~$X$ and~$Y$ of~${\mathfrak {ST}}_m$ will be said to be \textbf{adjacent} if and only if~$X$ and~$Y$ are two consecutive vertices of~${\mathfrak {ST}}_m$. We will write:

$$X \underset{m }{\sim}  Y$$

\noindent This edge relation ensures the existence of a word~{${\cal W}=\left ( {\cal W}_1, \hdots, {\cal W}_m\right)$} of length~$ m$, such that~$X$ and~$Y$ both belong to the iterate:

$$f_{\cal W} \,V_0=\left (f_{{\cal W}_1} \circ \hdots \circ  f_{{\cal W}_m} \right) \,V_0$$

 \noindent Given two points~$X$ and~$Y$ of the graph~$\mathfrak {ST}$, we will say that~$X$ and~$Y$ are \textbf{\emph{adjacent}} if and only if there exists a natural integer~$m$ such that:
$$X  \underset{m }{\sim}  Y$$
\end{definition}

 \vskip 1cm

\begin{proposition}\textbf{Adresses, on the the Sierpi\'{n}ski tetrahedron}\\

\noindent Given a strictly positive integer~$m$, and a word~${\cal W}=\left ( {\cal W}_1, \hdots, {\cal W}_m\right)$ of length~$m\,\in\,\N^\star$, on the graph~
${\mathfrak {ST}}_m  $, for any integer~$j$ of~$\left \lbrace0,...,3\right \rbrace  $, any~$X=f_{\cal W}(P_j)$ of~$ V_m \setminus V_{0}$, i.e. distinct from one of the~four fixed point~$P_i$, ~\mbox{$0 \leq i \leq 3$}, has exactly three adjacent vertices in $f_{\cal W}(V_0)$, given by:

$$f_{\cal W}(P_{{j+n}\pmod 4}) \quad \text{for} \ n\in\{1,2,3\}$$

\noindent where:

$$f_{\cal W}  = f_{{\cal W}_1} \circ \hdots \circ  f_{{\cal W}_m}   $$

\end{proposition}

\vskip 1cm

\begin{proposition}
Given a natural integer~$m$, a point~$X\, \in \,{\mathfrak {ST}}_m$ and a word~$\cal W $ of length~$m$ such that:

$$\mathcal{W}=\mathcal{W}_1\mathcal{W}_2...\mathcal{W}_{m-1}\mathcal{W}_m $$
$$X=F_{\mathcal{W}}(P_i) \, \in\,  V_m \setminus V_0 \quad , \quad i\,\in\, \{0,1,2,3\}$$

 \noindent let us write, for any integer~$\mathcal{W}_m\,\in \,\{0,1,2,3\}$n$\mathcal{W}_m:=j$, so:
$$\mathcal{W}=\mathcal{W}_1\mathcal{W}_2...\mathcal{W}_{m-1}j=\mathcal{\tilde{W}}j$$

Then $\mathcal{\tilde{W}}\in \{0,1,2,3\}^{m-1}$. The point~$X$ has exactly six adjacent vertices, of the form
$$f_{\mathcal{\tilde{W}}j}(P_k)$$
\noindent and
 $$f_{\mathcal{\tilde{W}}i}(P_l)$$
\noindent  for $k \neq i$ and $l \neq j$.\\
\end{proposition}

 \vskip 1cm
\begin{proof}
$x$ belongs to $f_{\mathcal{W}}(V_0)$, so it has three "neighbors" or adjacent vertices of the form $f_{\mathcal{\tilde{W}}j}(p_k)$ for $k\neq i$. And we recall that $f_{\mathcal{\tilde{W}}j}(p_i)=f_{\mathcal{\tilde{W}}i}(p_j)$, so $x$ belong to $f_{\mathcal{\tilde{W}}i}(V_0)$ too, so it has three other neighbors $f_{\mathcal{\tilde{W}}i}(p_l)$ for $l\neq i$.\\
\end{proof}

 \vskip 1cm

\begin{proposition}
Let us set:

$$\ds{V_\star =\underset{{m\in \N}}\bigcup \,V_m} $$

\noindent The set $V_{\star}$ is dense in~$\mathfrak{ST}$.\\
\end{proposition}

\section{Dirichlet forms on the Sierpi\'{n}ski tetrahedron}

Following J.~Kigami's approach \cite{Kigami1989}, Dirichlet forms and Laplacian on the Sierpi\'{n}ski tetrahedron can be respectively defined as limits of Dirichlet forms and Laplacians on~$\left (V_m\right)_{m \in\N}$.

\vskip 1cm

\begin{definition}\textbf{Dirichlet form, on a finite set} (\cite{Kigami2003})\\

 \noindent Let~$V$ denote a finite set~$V$, equipped with the usual inner product which, to any pair~$(u,v)$ of functions defined on~$V$, associates:

  $$(u,v)= \displaystyle \sum_{P\in  V} u(P)\,v(P)$$

  \noindent A \emph{\textbf{Dirichlet form}}on~$V$ is a symmetric bilinear form~${\cal E}$, such that:\\

\begin{enumerate}

\item For any real valued function~$u$ defined on~$V$:  ${\cal E}(u,u) \geq 0$.

\item   $  {\cal {E}}(u,u)= 0$ if and only if~$u$ is constant on~$V$.

\item For any real-valued function~$u$ defined on~$V$, if:
$$ u_\star = \min\, (\max(u, 0) , 1)  $$

\noindent i.e. :

$$\forall \,p \,\in\,V \, : \quad u_\star(p)= \left \lbrace \begin{array}{ccc} 1 & \text{if}& u(p) \geq 1 \\u(p) & \text{si}& 0 <u(p) < 1 \\0  & \text{if}& u(p) \leq 0 \end{array} \right.$$

\noindent then: ${ \cal{E}}(u_\star,u_\star)\leq { \cal{E}}(u,u)$ (Markov property).

\end{enumerate}

\end{definition}

\vskip 1cm

\begin{definition}\textbf{Energy, on the graph~${\mathfrak {ST}}_m $,~$m \,\in\,\N$, of a pair of functions}\\

 \noindent Let~$m$ be a natural integer, and~$u$ and~$v$ two real valued functions, defined on the set:

 $$V_m = \left \lbrace {\cal S}_0^m,  {\cal S}_1^m, \hdots,  {\cal S}_{{\cal N}_m -1}^m \right \rbrace $$

 \noindent of the~${\cal N}_m $ vertices of~${\mathfrak {ST}}_m$.\\

 \noindent \textbf{The energy, on the graph~${\mathfrak {ST}}_m$, of the pair of functions~$(u,v)$, is:}

\[ E_{{\mathfrak {ST}}_m}(u,v)=\displaystyle \sum_{i=0}^{\mathcal{N}_m - 2}\left( u(\mathcal{S}_i^m)-u(\mathcal{S}_{i+1}^m)\right) \left( v(\mathcal{S}_i^m)-v(\mathcal{S}_{i+1}^m)\right)  \]
or
\[ E_{{\mathfrak {ST}}_m}(u,v)=\displaystyle\sum_{X \underset{m}{\sim} Y}\left( u(X)-u(Y)\right) \left( v(X)-v(Y)\right)  \]
\noindent Let us note that
$$E_{{\mathfrak {ST}}_m}(u,u)=0 \quad \text{ if $u$ is constant}$$

\noindent  $E_{{\mathfrak {ST}}_m}$ is a Dirichlet form on ${\mathfrak {ST}}_m$.\\
\end{definition}

\vskip 1cm

\begin{proposition}\textbf{Harmonic extension of a function, on the Sierpi\'{n}ski Tetrahedron}\\

\noindent For any strictly positive integer~$m$, if~$u$ is a real-valued function defined on~$V_{m-1}$, its \textbf{harmonic extension}, denoted by~$ \tilde{u}$, is obtained as the extension of~$u$ to~$V_m$ which minimizes the energy:

$$  {\cal{E}}_{{\mathfrak {ST}}_m }(\tilde{u},\tilde{u})=\displaystyle \sum_{X \underset{m }{\sim} Y} (\tilde{u}(X)-\tilde{u}(Y))^2 $$

\noindent The link between~$   {\cal{E}}_{{\mathfrak {ST}}_m }$ and~$  {\cal{E}}_{{\mathfrak {ST}}_{m-1}}$ is obtained through the introduction of two strictly positive constants~$r_{ m}$,~$r_{m+1}$, such that:

$$  r_{m }\, \displaystyle \sum_{X \underset{m  }{\sim} Y} (\tilde{u}(X)-\tilde{u}(Y))^2 = r_{m-1}\,\displaystyle  \sum_{X \underset{m-1 }{\sim} Y} (u(X)-u(Y))^2$$

\noindent In particular:

$$  r_{1 }\, \displaystyle \sum_{X \underset{1  }{\sim} Y} (\tilde{u}(X)-\tilde{u}(Y))^2 = r_{0}\,\displaystyle  \sum_{X \underset{0 }{\sim} Y} (u(X)-u(Y))^2$$

\noindent For the sake of simplicity, we will fix the value of the initial constant:~$r_0=1$. One has then:

$$ {\cal{E}}_{ {\mathfrak {ST}}_m} (\tilde{u},\tilde{u})= \displaystyle \frac{1}{ r_{1 }}\,  {\cal{E}}_{{\mathfrak {ST}}_0}(\tilde{u},\tilde{u})$$

\noindent Let us set:

$$r = \displaystyle \frac{1}{r_{1 }} $$

\noindent and:

$$  {\cal{E}}_{m}(u)=r_m\, \sum_{X \underset{m }{\sim} Y} (\tilde{u}(X)-\tilde{u}(Y))^2 $$

\noindent Since the determination of the harmonic extension of a function appears to be a local problem, on the graph~${\mathfrak {ST}}_{m-1}$, which is linked to the graph~${\mathfrak {ST}}_m$ by a similar process as the one that links~${\mathfrak {ST}}_1$ to~${\mathfrak {ST}}_0$, one deduces, for any strictly positive integer~$m$:

$$ {\cal{E}}_{{\mathfrak {ST}}_m}(\tilde{u},\tilde{u})= \displaystyle \frac{1}{ r_{1 }}\,  {\cal{E}}_{{\mathfrak {ST}}_{m-1}}(\tilde{u},\tilde{u})$$

\noindent By induction, one gets:

$$r_m=r_1^m\,r_0=r^{-m} $$

\noindent If~$v$ is a real-valued function, defined on~$V_{m-1}$, of harmonic extension~$ \tilde{v}$, we will write:

$$  {\cal{E}}_{m}(u,v)=r^{-m}\, \sum_{X \underset{m }{\sim} Y} (\tilde{u}(X)-\tilde{u}(Y)) \, (\tilde{v}(X)-\tilde{v}(Y)) $$

\noindent For further precision on the construction and existence of harmonic extensions, we refer to~\emph{\cite{Sabot1987}}.
\end{proposition}

\vskip 1cm

\begin{definition}\textbf{Dirichlet form, for a pair of continuous functions defined on the Sierpi\'{n}ski tetrahedron~$\mathfrak {ST}$}\\

 \noindent We define the Dirichlet form~$\cal{E}$ which, to any pair of real-valued, continuous functions~$(u,v)$ defined on~$\mathfrak {ST}$, associates, subject to its existence:

$$
  {\cal{E}} (u,v)= \displaystyle \lim_{m \to + \infty} {\cal{E}}_{m }\left (u_{\mid V_m},v_{\mid V_m}\right)=
  \displaystyle \lim_{m \to + \infty}\displaystyle \sum_{X  \underset{m }{\sim}  Y} r^{-m}\,  \left (u_{\mid V_m}(X)-u_{\mid V_m}(Y)\right )\,\left(v_{\mid V_m}(X)-v_{\mid V_m}(Y)\right) $$

\end{definition}

\vskip 1cm
\begin{definition}\textbf{Normalized energy, for a continuous function~$u$, defined on~$\mathfrak {ST} $}\\
\noindent Taking into account that the sequence~$\left (\mathcal{E}_m\left ( u_{\mid V_m} \right)\right)_{m\in\N}$ is defined on
$$\ds{V_\star =\underset{{i\in \N}}\bigcup \,V_i}$$

\noindent one defines the normalized energy, for a continuous function~$u$, defined on~$\mathfrak {ST}$, by:

$$\mathcal{E}(u)=\underset{m\rightarrow +\infty}\lim \mathcal{E}_m \left ( u_{\mid V_m} \right)$$

\end{definition}

\vskip 1cm

\begin{pte}
 \noindent The Dirichlet form~$\cal{E}$ which, to any pair of real-valued, continuous functions defined on~$\mathfrak {ST}$, associates:

$$
  {\cal{E}} (u,v)= \displaystyle \lim_{m \to + \infty} {\cal{E}}_{m }\left (u_{\mid V_m},v_{\mid V_m}\right)=
  \displaystyle \lim_{m \to + \infty}\displaystyle \sum_{X  \underset{m }{\sim}  Y} r^{-m}\,\left (u_{\mid V_m}(X)-u_{\mid V_m}(Y)\right )\,\left(v_{\mid V_m}(X)-v_{\mid V_m}(Y)\right) $$

\noindent satisfies the self-similarity relation:

$$
  {\cal{E}} (u,v)= r^{-1}\,\displaystyle \sum_{i=0}^{3}   {\cal{E}} \left ( u \circ f_i,v\circ f_i \right)$$

\end{pte}

\vskip 1cm

\begin{proof}

$$\begin{array}{ccc}  \displaystyle \sum_{i=0}^{3}   {\cal{E}} \left ( u \circ f_i,v\circ f_i \right)&=& \displaystyle \lim_{m \to + \infty} \displaystyle \sum_{i=0}^{3}  {\cal{E}}_{m  }\left ( u_{\mid V_{m }} \circ f_i,v_{\mid V_{m }}\circ f_i \right)\\
&=&\displaystyle \lim_{m \to + \infty}\displaystyle \sum_{X  \underset{m  }{\sim}  Y} r^{-m }\, \displaystyle \sum_{i=0}^{3}   \left (u_{\mid V_{m }}\left (f_i(X)\right )-
u_{\mid V_{m }} \left (f_i(Y)\right  )\right )\,\left(v_{\mid V_{m }}\left (f_i(X)\right )-v_{\mid V_{m }}\left (f_i(Y)\right )\right) \\
&=&\displaystyle \lim_{m \to + \infty}\displaystyle \sum_{X  \underset{m+1  }{\sim}  Y} r^{-m }\, \displaystyle \sum_{i=0}^{3}   \left (u_{\mid V_{m+1}} (X) -u_{\mid V_{m+1}}  (Y) \right )\,
\left(v_{\mid V_{m+1}} (X) -v_{\mid V_{m+1}} (Y) \right) \\
&=&\displaystyle \lim_{m \to + \infty}r \,{\cal{E}}_{m+1  }\left (u_{\mid V_{m+1}},v_{\mid V_{m+1}} \right) \\
&=& r\, {\cal{E}} (u,v) \\
\end{array} $$

\end{proof}
\vskip 1cm

\begin{notation}
\noindent We will denote by~$\text{dom}\,{\cal E}$ the subspace of continuous functions defined on~$\mathfrak {ST}$, such that:

$$\mathcal{E}(u)< + \infty$$

\end{notation}

\vskip 1cm

\begin{notation}
\noindent We will denote by~$\text{dom}_0\,{ \cal E}$ the subspace of continuous functions defined on~$\mathfrak {ST}$, which take the value $0$ on~$V_0$, such that:

$$\mathcal{E}(u)< + \infty$$

\end{notation}

\vskip 1cm
\begin{proposition}
The space $\text{dom}\mathcal{E}$, modulo the space of constant function on $\mathfrak{ST}$, is a Hilbert space.\\
\end{proposition}

\section{Explicit construction of the Dirichlet forms}
  Let us denote by~$u$ a real valued function defined on:
$$V_0=\{P_0,P_1,P_2,P_3\}$$

\noindent We herafter aim at determining its harmonic extension~$\tilde{u}$ on~$V_1$.\\

\noindent For the sake of simplicity, we set:

$$u(p_0)=a \quad , \quad u(p_1)=b \quad , \quad u(p_2)=c\quad , \quad u(p_3)=d$$

\noindent One has to bear in mind that the energy on~$V_0$ is given by:
\[E_0(u)=(a-b)^2+(a-c)^2+(a-d)^2+(b-c)^2+(b-d)^2+(c-d)^2\]

\noindent For the sake of simplicity, we set:
$$\tilde{u}(f_0(q_1))=x_1 \quad , \quad \tilde{u}(f_1(q_2))=x_2 \quad , \quad\tilde{u}(f_0(q_2))=x_3 \quad , \quad\tilde{u}(f_0(q_3))=x_4\quad , \quad \tilde{u}(f_1(q_3))=x_5\quad , \quad\tilde{u}(f_2(q_3))=x_6$$

\noindent Thus:

\begin{align*}
E_1(\tilde{u})&=(x_1-a)^2+(x_1-b)^2+(x_1-x_2)^2+(x_1-x_3)^2+(x_1-x_4)^2+(x_1-x_5)^2\\
&+(x_2-b)^2+(x_2-c)^2+(x_2-x_3)^2+(x_2-x_5)^2+(x_2-x_6)^2\\
&+(x_3-a)^2+(x_3-c)^2+(x_3-x_4)^2+(x_3-x_6)^2\\
&+(x_4-a)^2+(x_4-d)^2+(x_4-x_5)^2+(x_4-x_6)^2\\
&+(x_5-b)^2+(x_5-d)^2+(x_5-x_6)^2\\
&+(x_6-c)^2+(x_6-d)^2
\end{align*}
\noindent The minimum of this quantity is to be obtained in the set of critical points, which leads to:

\begin{gather*}
6x_1-x_2-x_3-x_4-x_5=a+b \\
6x_2-x_1-x_3-x_5-x_6=b+c \\
6x_3-x_1-x_2-x_4-x_6=a+c \\
6x_4-x_1-x_3-x_5+x_6=a+d \\
6x_5-x_1-x_2-x_4-x_6=b+d \\
6x_6-x_2-x_3-x_4-x_5=c+d \\
\end{gather*}
\noindent Under matricial form:
\[ \mathbf{x}=\mathbf{A}^{-1}\mathbf{b}\]
\noindent with
\[
\mathbf{A} =
\begin{pmatrix}
6 & -1 & -1 & -1 & -1 & 0 \\
-1 & 6 & -1 & 0 & -1 & -1 \\
-1 & -1 & 6 & -1 & 0 & -1 \\
-1 & 0 & -1 & 6 & -1 & -1 \\
-1 & -1 & 0 & -1 & 6 & -1 \\
0 & -1 & -1 & -1 & -1 & 6 \\
\end{pmatrix}\]

\[
\mathbf{x} =
\begin{pmatrix}
x_1 \\
x_2 \\
x_3 \\
x_4 \\
x_5 \\
x_6 \\
\end{pmatrix}\]

\[
\mathbf{b} =
\begin{pmatrix}
a+b \\
b+c \\
a+c \\
a+d \\
b+d \\
c+d \\
\end{pmatrix}\]

\noindent Finally, we get:

\[
\mathbf{x}=
\left(
\begin{array}{c}
 \frac{1}{6} (2 a+2 b+c+d) \\
 \frac{1}{6} (a+2 b+2 c+d) \\
 \frac{1}{6} (2 a+b+2 c+d) \\
 \frac{1}{6} (2 a+b+c+2 d) \\
 \frac{1}{6} (a+2 b+c+2 d) \\
 \frac{1}{6} (a+b+2 (c+d)) \\
\end{array}
\right)
\]

\noindent By substituting~$\mathbf{X}$ in the energy, one obtains:
\begin{align*}
E_1(\tilde{u}) &=\left\{\displaystyle \frac{2}{3} \left(3 a^2-2 a (b+c+d)+3 b^2-2 b (c+d)+3 c^2-2 c d+3
   d^2\right)\right\} \\
&=\displaystyle\frac{2}{3}\, E_0(u)
\end{align*}

\noindent Let us now move to the general case, and consider a natural integer~$m  $. Each point of~\mbox{$V_{m+1} \setminus V_{m}$} belongs to a~$m$-cell of the form $$f_{\mathcal{W}}(\mathfrak{ST})$$

\noindent where~$ \mathcal{W} $ denotes a word of length~$m$. The total energy~$E_{m+1}(\tilde{u})$ is given by:
\[E_{m+1}(\tilde{u})=\sum_{\left|\mathcal{W}\right|=m}E_1(\tilde{u}\circ f_{\mathcal{W}})\]
\noindent The global minimization problem can be reduced to~$4^m$ local minimization problems which are of the same kind of the one we just solved.\\
\noindent Thus, the normalization constant is.

$$\displaystyle{r=\displaystyle \frac{2}{3}}$$

\noindent This enables us to define the normalized energy:
\[ \mathcal{E}_m(u)=r^{-m}\,E_m(u) \]
\noindent and its limit:
\[\mathcal{E}(u)=\displaystyle \lim_{m\rightarrow +\infty}\mathcal{E}_m(u) \]
for $u\,\in \,\text{dom}(\mathcal{E})$.\\

 \begin{figure}[h!]
 \center{\psfig{height=8cm,width=10cm,angle=0,file=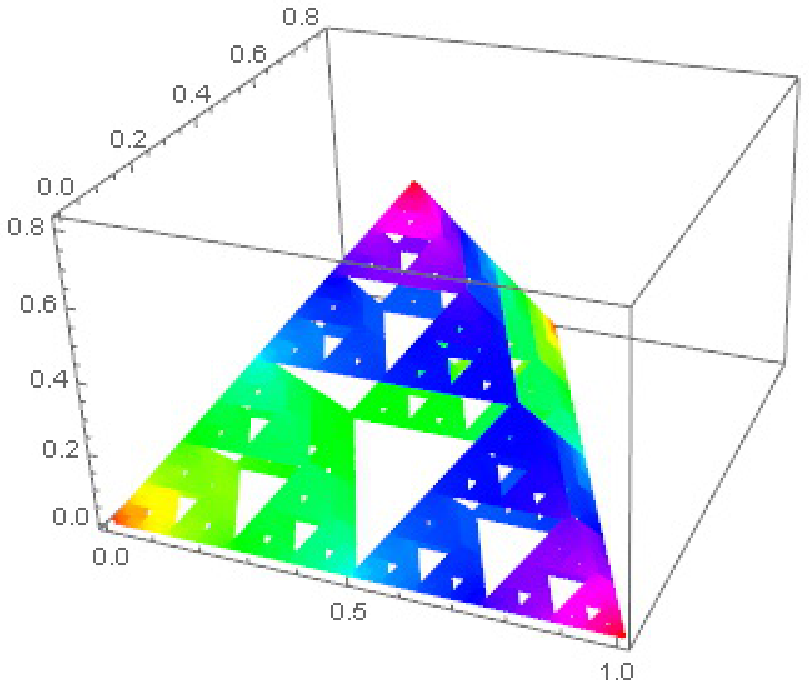} \psfig{height=8cm,width=1cm,angle=0,file=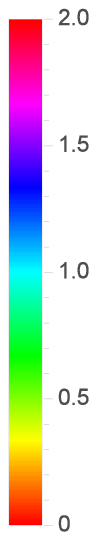}}
\caption{The harmonic extension on the Sierpi\'{n}ski tetrahedron of a function taking the values~\mbox{$a=0$,~$b=2$,~$c=0$ and~$d=2$}.}
\end{figure}

\newpage
\section{Laplacian, on the Sierpi\'{n}ski tetrahedron}

\begin{definition}\textbf{Self-similar measure, on the Sierpi\'{n}ski tetrahedron}\\

\noindent A measure~$\mu$ on~$\R^3$ will be said to be \textbf{self-similar} on the Sierpi\'{n}ski tetrahedron, if there exists a family of strictly positive pounds~\mbox{$\left (\mu_i\right)_{0 \leq i \leq 3}$} such that:

$$ \mu= \displaystyle \sum_{i=0}^{3} \mu_i\,\mu\circ f_i^{-1} \quad, \quad \displaystyle \sum_{i=0}^{3} \mu_i =1$$

\noindent For further precisions on self-similar measures, we refer to the works of~J.~E.~Hutchinson~(see \cite{Hutchinson1981}).

\end{definition}

\vskip 1cm

\begin{pte}\textbf{Building of a self-similar measure, for the Sierpi\'{n}ski tetrahedron}\\

\noindent The Dirichlet forms mentioned in the above require a positive Radon measure with full support.

\noindent Let us set, for any integer~$i$ belonging to~$\left \lbrace 0, \hdots, 3 \right \rbrace$:

 $$\mu_i=\displaystyle \frac{1}{4}$$

  \noindent  This enables one to define a self-similar measure~$\mu$ on~$\mathfrak{ST}$ as:
\[\mu =\displaystyle \frac{1}{4}\,\sum_{i=0}^3  \mu \circ f_i \]

\end{pte}

\vskip 1cm

\begin{definition}\textbf{Laplacian of order~$m\,\in\,\N^\star$}\\

\noindent For any strictly positive integer~$m$, and any real-valued function~$u$, defined on the set~$V_m$ of the vertices of the graph~${\mathfrak {ST}}_m $, we introduce the Laplacian of order~$m$,~$\Delta_m(u)$, by:

$$\Delta_m u(X) = \displaystyle\sum_{Y \in V_m,\,Y\underset{m}{\sim} X} \left (u(Y)-u(X)\right)  \quad \forall\, X\,\in\, V_m\setminus V_0 $$

\end{definition}

\vskip 1cm

\begin{definition}\textbf{Harmonic function of order~$m\,\in\,\N^\star$}\\

\noindent Let~$m$ be a strictly positive integer. A real-valued function~$u$,defined on the set~$V_m$ of the vertices of the graph~${\mathfrak {ST}}_m $, will be said to be \textbf{harmonic of order~$m$} if its Laplacian of order~$m$ is null:

$$\Delta_m u(X) =0 \quad \forall\, X\,\in\, V_m\setminus V_0 $$

\end{definition}

\vskip 1cm

\begin{definition}\textbf{Piecewise harmonic function of order~$m\,\in\,\N^\star$}\\

\noindent  Given a strictly positive integer~$m$, a real valued function~$u$, defined on the set of vertices of~$\mathfrak {ST}$, is said to be \textbf{piecewise harmonic function of order~$m$} if, for any word~${\cal W}$ of length~$ m$,~$u\circ f_{\cal W}$ is harmonic of order~$m$.

\end{definition}

\vskip 1cm

\begin{definition}\textbf{Existence domain of the Laplacian, for a continuous function on~$\mathfrak {ST}$} (see \cite{Beurling1985})\\

\label{Lapl}
\noindent We will denote by~$\text{dom}\, \Delta$ the existence domain of the Laplacian, on the graph~${\mathfrak {ST}}$, as the set of functions~$u$ of~$\text{dom}\, \mathcal{E}$such that there exists a continuous function on~$\mathfrak {ST}$, denoted~$\Delta \,u$, that we will call \textbf{Laplacian of~$u$}, such that :
$$\mathcal{E}(u,v)=-\displaystyle \int_{{\cal D} \left ( {\mathfrak {ST} }\right)} v\, \Delta u   \,d\mu \quad \text{for any } v \,\in \,\text{dom}_0\, \mathcal{E}$$
\end{definition}

\vskip 1cm

\begin{definition}\textbf{Harmonic function}\\

\noindent A function~$u$ belonging to~\mbox{$\text{dom}\,\Delta$} will be said to be \textbf{harmonic} if its Laplacian is equal to zero.
\end{definition}

\vskip 1cm

\begin{notation}

In the following, we will denote by~${\cal H}_0\subset \text{dom}\, \Delta$ the space of harmonic functions, i.e. the space of functions~$u \,\in\,\ \text{dom}\, \Delta$ such that:

$$\Delta\,u=0$$

\noindent Given a natural integer~$m$, we will denote by~${\cal S} \left ({\cal H}_0,V_m \right)$ the space, of dimension~$4^m$, of spline functions " of level~$m$", ~$u$, defined on~${\mathfrak {ST}}$, continuous, such that, for any word~$\cal W$ of length~$m$,~\mbox{$u \circ T_{\cal W}$} is harmonic, i.e.:

$$\Delta_m \, \left ( u \circ T_{\cal W} \right)=0$$

\end{notation}

\vskip 1cm

\begin{pte}

For any natural integer~$m$:

$${\cal S} \left ({\cal H}_0,V_m \right )\subset  \text{dom }{ \cal E}$$

\end{pte}
\vskip 1cm

\begin{pte}
Let~$m$ be a strictly positive integer,~$X \,\notin\,V_0$ a vertex of the graph~$\mathfrak {ST}$, and~\mbox{$\psi_X^{m}\,\in\,{\cal S} \left ({\cal H}_0,V_m \right)$} a spline  function such that:

$$\psi_X^{m}(Y)=\left \lbrace \begin{array}{ccc}\delta_{XY} & \forall& Y\,\in \,V_m \\
 0 & \forall& Y\,\notin \,V_m \end{array} \right. \quad,  \quad \text{where} \quad    \delta_{XY} =\left \lbrace \begin{array}{ccc}1& \text{if} & X=Y\\ 0& \text{else} &  \end{array} \right.$$

\noindent Then, since~$X\, \notin \,V_0$: $\psi_X^{m}\,\in \,\text{dom}_0\, \mathcal{E}$.

\noindent For any function~$u$ of~$\text{dom}\, \mathcal{E}$, such that its Laplacian exists, definition (\ref{Lapl}) applied to~$\psi_X^{m}$ leads to:

$$\mathcal{E}(u,\psi_X^{m})=\mathcal{E}_m(u,\psi_X^{m})= -r^{-m}\,\Delta_m u(X)=- \displaystyle\int_{{\cal D}\left ({\mathfrak {ST}}\right)}  \psi_X^{m}\,\Delta u  \, d\mu \approx -\Delta  u(X)\, \displaystyle\int_{{\cal D} \left ( {\mathfrak {ST}} \right)}  \psi_X^{m}\, d\mu$$

\noindent since~$\Delta u$ is continuous on~$ {\mathfrak {ST}}$, and the support of the spline function~$\psi_X^{m}$ is close to~$X$:

$$\displaystyle\int_{{\cal D} \left ( {\mathfrak {ST}}  \right)}  \psi_X^{m}\,\Delta u  \, d\mu \approx -\Delta  u(X)\, \displaystyle\int_{{\cal D} \left ( {\mathfrak {ST}} \right)}  \psi_X^{m}\, d\mu$$

\noindent By passing through the limit when the integer~$m$ tends towards infinity, one gets:

$$ \displaystyle \lim_{m \to + \infty} \displaystyle\int_{{\cal D} \left ( {\mathfrak {ST}} \right)}  \psi_X^{m}\,\Delta_m u  \, d\mu=
 \Delta  u(X)\,\displaystyle \lim_{m \to + \infty}   \displaystyle\int_{{\cal D} \left ( {\mathfrak {ST}} \right)}  \psi_X^{m}\, d\mu$$

\noindent i.e.:

$$  \Delta  u(X)= \displaystyle \lim_{m \to + \infty} r^{-m} \left (  \displaystyle\int_{{\cal D} \left ( {\mathfrak {ST}}  \right)}  \psi_X^{m}\, d\mu  \right)^{-1} \,\Delta_m u(X)\,$$

\end{pte}

\vskip 1cm

\begin{remark}

\noindent As it is explained in~\cite{StrichartzLivre2006}, one has just to reason by analogy with the dimension~1, more particulary, the unit interval~$I=[0,1]$, of extremities~$X_0=(0,0)$, and~$X_1=(1,0)$. The functions~$\psi_{X_1}$ and~$\psi_{X_2}$ such that, for any~$Y$ of~$\R^2$ :

$$\psi_{X_1} (Y)=\delta_{X_1Y} \quad  ,  \quad \psi_{X_2} (Y)=\delta_{X_2Y}   $$

\noindent are, in the most simple way, tent functions. For the standard measure, one gets values that do not depend on~$X_1$, or~$X_2$ (one could, also, choose to fix~$X_1$ and~$X_2$ in the interior of~$I$) :

$$\displaystyle\int_{I}  \psi_{X_1}\, d\mu =\displaystyle\int_{I}  \psi_{X_2}\, d\mu=\displaystyle \frac{1}{2}$$

\noindent (which corresponds to the surfaces of the two tent triangles.) \\

 \begin{figure}[h!]
 \center{\psfig{height=8cm,width=10cm,angle=0,file=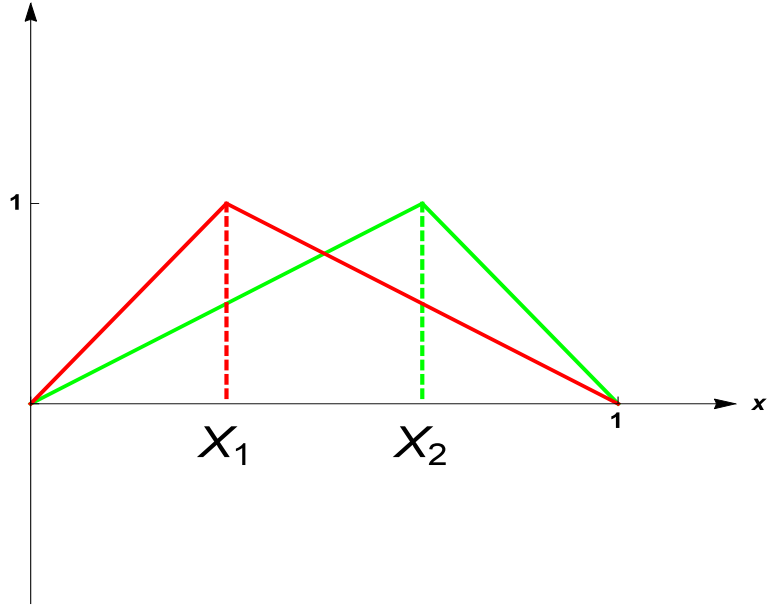}}\\
\caption{The graphs of the spline functions~$\psi_{X_1}$ and~$\psi_{X_2}$.}
 \end{figure}

\noindent In our case, we have to build the pendant, we no longer reason on the unit interval, but on our polyhedra cells. \\

\noindent Given a natural integer~$m$, and a point~$X \,\in\,V_m$, the spline function~$\psi_X^{m}$ is supported by two~$m$-polyhedra cells. It is such that, for every $m$-polyhedra cell $f_{\mathcal{W}} \left( \mathfrak{ST}\right)$ the vertices of which are~$X$,~\mbox{$Y\neq X$},\\~\mbox{$Z\neq X$},~\mbox{$T\neq X$}:
$$\psi_X^{m}+\psi_Y^{m}+\psi_Z^{m}+\psi_t^{m}=1$$

\noindent Thus:
\[ \displaystyle \int_{f_{\mathcal{W}} \left( \mathfrak{ST}\right)} \left(\psi_X^{m}+\psi_Y^{m}+\psi_Z^{m}+\psi_T^{m} \right)\, d\mu=\mu(f_{\mathcal{W}} \left( \mathfrak{ST}\right) )=\displaystyle \frac{1}{4^m} \]
\noindent By symmetry, all three summands have the same integral. This yields:

$$\displaystyle{\int_{f_{\mathcal{W}} \left( \mathfrak{ST}\right)} \psi_X^{m} d\mu=\frac{1}{4^{m+1}}}$$
\noindent Taking into account the contributions of the remaining~$m$-polyhedra cells, one has:

$$\displaystyle{\int_{\mathfrak{ST}} \psi_X^{m} d\mu=\frac{2}{4^{m+1}}}$$

\noindent which leads to:

$$\displaystyle{\left (\int_{\mathfrak{ST}} \psi_X^{m}\, d\mu \right )^{-1}=\displaystyle\frac{4^{m+1}}{2}} $$

\noindent Since:

$$ \displaystyle{r^{-m}=\left(\displaystyle \frac{3}{2} \right)^m}$$

\noindent this enables us to obtain the point-wise formula, for~$u\,\in\,\text{dom}\,\Delta$:

$$\forall \,X \,\in\,{\mathfrak {ST}} \,: \quad  \Delta_{\mu} u(X)=2 \displaystyle\lim_{m\rightarrow +\infty} 6^m\,
\Delta_m u(X) $$

\end{remark}

\vskip 1cm

\vskip 1cm
\begin{theorem}

\noindent Let~$u$ be in~\mbox{$\text{dom}\,\Delta$}. Then, the sequence of functions~$\left (f_m \right)_{m\in\N}$ such that, for any natural integer~$m$, and any~$X$ of~\mbox{$V_\star\setminus V_0$} :

 $$f_m(X)=r^{-m}\,\left (\displaystyle \int_{{\cal D} \left ({\mathfrak {ST}}\right) }   \psi_{X}^{m}\,d\mu\right)^{-1}\,\Delta_m \,u(X) $$

 \noindent  converges uniformly towards~$\Delta\,u$, and, reciprocally, if the sequence of functions~$\left (f_m \right)_{m\in\N}$ converges uniformly towards a continuous function on~\mbox{$V_\star\setminus V_0$}, then:

 $$u \,\in\, \text{dom}\,\Delta$$
\end{theorem}

\vskip 1cm
\begin{proof}

\noindent Let~$u$ be in~\mbox{$\text{dom}\,\Delta$}. Then:

 $$ r^{-m}\,\left (\displaystyle \int_{{\cal D} \left ({\mathfrak {ST}}\right) }  \psi_{X}^{m}\,d\mu\right)^{-1}\,\Delta_m \,u(X)=
\displaystyle \frac{\displaystyle \int_{{\cal D} \left ({\mathfrak {ST}}\right) }  \Delta\,u \,\psi_{X}^{m}\,d\mu} {\displaystyle \int_{{\cal D} \left ({\mathfrak {ST}}\right) }   \psi_{X}^{m}\,d\mu}   $$

\noindent Since~$u$ belongs to~\mbox{$\text{dom}\,\Delta$}, its Laplacian~$\Delta\,u$ exists, and is continuous on the graph~${\mathfrak {ST}}$. The uniform convergence of the sequence~$\left (f_m \right)_{m\in\N}$ follows.\\

\noindent Reciprocally, if the sequence of functions~$\left (f_m \right)_{m\in\N}$ converges uniformly towards a continuous function on~\mbox{$V_\star\setminus V_0$}, the, for any natural integer~$m$, and any~$v$ belonging to~\mbox{$\text{dom}_0\,{\cal E}$}:

$$\begin{array}{ccc} {\cal{E}}_{m }(u,v)
  &=&  \displaystyle \sum_{(X,Y) \,\in \, V_m^2,\,X  \underset{m }{\sim}  Y} r^{-m}\,\left (u_{\mid V_m}(X)-u_{\mid V_m}(Y)\right )\,\left(v_{\mid V_m}(X)-v_{\mid V_m}(Y)\right) \\
  &=& \displaystyle \sum_{(X,Y) \,\in \, V_m^2,\,X  \underset{m }{\sim}  Y} r^{-m}\,\left (u_{\mid V_m}(Y)-u_{\mid V_m}(X )\right )\,\left(v_{\mid V_m}(Y)-v_{\mid V_m}(X)\right) \\
  &=&- \displaystyle \sum_{X \,\in \,V_m\setminus V_0 } r^{-m}\,\sum_{Y\,\in \,V_m, \, Y  \underset{m }{\sim}  X} v_{\mid V_m}(X)\,\left (u_{\mid V_m}(Y)-u_{\mid V_m}(X)\right )   \\
  & &- \displaystyle \sum_{X \,\in \,  V_0 } r^{-m}\,\sum_{Y\,\in \,V_m ,\, Y  \underset{m }{\sim}  X} v_{\mid V_m}(X)\,\left (u_{\mid V_m}(Y)-u_{\mid V_m}(X)\right )   \\
  &=&- \displaystyle \sum_{X \,\in \,V_m\setminus V_0 } r^{-m}\,v(X)\,\Delta_m \,u(X)  \\
    &=&- \displaystyle \sum_{X \,\in \,V_m\setminus V_0 } v(X)\,   \left (\displaystyle \int_{{\cal D} \left ({\mathfrak {ST}}\right) }   \psi_{X}^{m}\,d\mu\right) \, r^{-m}\, \left (\displaystyle \int_{{\cal D} \left ({\mathfrak {ST}}\right) }   \psi_{X}^{m}\,d\mu\right)^{-1}\, \Delta_m \,u(X)  \\
  \end{array}
  $$

\noindent Let us note that any~$X$ of~$V_m\setminus V_0$ admits exactly three adjacent vertices which belong to~$V_m\setminus V_0$, which accounts for the fact that the sum

 $$\displaystyle \sum_{X \,\in \,V_m\setminus V_0 } r^{-m}\,\sum_{Y\,\in \,V_m\setminus V_0 ,\, Y  \underset{m }{\sim}  X} v(X)\,\left (u_{\mid V_m}(Y)-u_{\mid V_m}(X)\right)$$
 \noindent has the same number of terms as:

 $$ \displaystyle \sum_{(X,Y) \,\in \,(V_m\setminus V_0)^2,\,X  \underset{m }{\sim}  Y} r^{-m}\,\left (u_{\mid V_m}(Y)-u_{\mid V_m}(X)\right )\,\left(v_{\mid V_m}(Y)-v_{\mid V_m}(X)\right)  $$

 \noindent For any natural integer~$m$, we introduce the sequence of functions~$\left (f_m \right)_{m\in\N}$ such that, for any~$X$ of~$V_m\setminus V_0$:

 $$f_m(X)=r^{-m}\,\left (\displaystyle \int_{{\cal D} \left ({\mathfrak {ST}} \right) }  \psi_{X}^{m}\,d\mu\right)^{-1}\,\Delta_m \,u(X) $$

 \noindent The sequence~$\left (f_m \right)_{m\in\N}$ converges uniformly towards~$\Delta\,u$. Thus:

$$\begin{array}{ccc} {\cal{E}}_{m }(u,v)
  &=&      -   \displaystyle \int_{{\cal D} \left ({\mathfrak {ST}}\right) }  \left \lbrace \displaystyle \sum_{X \,\in \,V_m\setminus V_0 } v_{\mid V_m}(X)\,\Delta\,u_{\mid V_m}(X)\, \psi_{X}^{m}\right \rbrace \,d\mu
  \end{array}
  $$

\end{proof}

\section{Normal derivatives}

Let us go back to the case of a function~$u$ twice differentiable on~$I=[0,1]$, that does not vanish in~0 and~:
$$  \displaystyle \int_0^1 \left (\Delta u \right)(x)\,v(x)\,dx=-   \displaystyle \int_0^1 u'(x)\,v'(x)\,dx+ u '(1)\,v (1)-u '(0)\,v (0) $$

\noindent The normal derivatives:

$$\partial_n u(1)=u'(1) \quad \text{et} \quad \partial_n u(0)=u'(0) $$

\noindent appear in a natural way. This leads to:

$$  \displaystyle \int_0^1 \left (\Delta u \right)(x)\,v(x)\,dx=-   \displaystyle \int_0^1 u'(x)\,v'(x)\,dx+ \displaystyle \sum_{\partial\, [0,1]} v\,\partial_n\,u $$

\noindent One meets thus a particular case of the Gauss-Green formula, for an open set~$\Omega$ of~$\R^d$,~$d \,\in\,\N^\star$:

$$  \displaystyle \int_\Omega \nabla\,  u\, \nabla \, v \, d \mu= -\displaystyle \int_\Omega  \left (\Delta u \right) \,v \,d\mu + \displaystyle \int_{\partial\, \Omega } v\,\partial_n\,u \,d\sigma$$

\noindent where~$\mu$ is a measure on~$\Omega $, and where~$d\sigma$ denotes the elementary surface on~$\partial\, \Omega $.\\

\noindent In order to obtain an equivalent formulation in the case of the graph~$\mathfrak {ST} $, one should have, for a pair of functions~$(u,v)$ continuous on~$\mathfrak {ST} $ such that~$u$ has a normal derivative:

$$  {\cal E}(u,v)= -\displaystyle \int_\Omega  \left (\Delta u \right) \,v \,d\mu + \displaystyle \sum_{V_0} v\,\partial_n\,u  $$

\noindent For any natural integer~$m$ :

$$\begin{array}{ccc}
 {\cal{E}}_{m }(u,v)\\
  &=&  \displaystyle \sum_{(X,Y) \,\in \, V_m^2,\,X  \underset{m }{\sim}  Y} r^{-m}\,\left (u_{\mid V_m}(Y)-u_{\mid V_m}(X)\right )\,\left(v_{\mid V_m}(Y)-v_{\mid V_m}(X)\right) \\
  &=&- \displaystyle \sum_{X \,\in \,V_m\setminus V_0 } r^{-m}\,\sum_{Y\,\in \,V_m  ,\, Y  \underset{m }{\sim}  X} v_{\mid V_m}(X)\,\left (u_{\mid V_m}(Y)-u_{\mid V_m}(X)\right )   \\
   &&-  \displaystyle \sum_{X \,\in \,  V_0 } r^{-m}\,\sum_{Y\,\in \,  V_m ,\, Y  \underset{m }{\sim}  X} v_{\mid V_m}(X)\,\left (u_{\mid V_m}(Y)-u_{\mid V_m}(X)\right )  \\
  &=&- \displaystyle \sum_{X \,\in \,V_m\setminus V_0 } v_{\mid V_m}(X)\,r^{-m}\,\Delta_m \,u_{\mid V_m}(X)  \\
   &&+ \displaystyle \sum_{X \,\in \,  V_0 } \sum_{Y\,\in \,  V_m ,\, Y  \underset{m }{\sim}  X} r^{-m}\,v_{\mid V_m}(X)\,\left (u_{\mid V_m}(X)-u_{\mid V_m}(Y)\right )  \\
  \end{array}
  $$

 \noindent We thus come across an analogous formula of the Gauss-Green one, where the role of the normal derivative is played by:

 $$ \displaystyle \sum_{X \,\in \,  V_0 } r^{-m}\,\sum_{Y\,\in \,  V_m ,\, Y  \underset{m }{\sim}  X} \,\left (u_{\mid V_m}(X)-u_{\mid V_m}(Y)\right ) $$

 \vskip 1cm

 \begin{definition}

 \noindent For any~$X$ of~$V_0$, and any continuous function~$u$ on~$\mathfrak {ST} $, we will say that~$u$ admits a normal derivative in~$X$, denoted by~$\partial_n\,u(X)$, if:

 $$ \displaystyle \lim_{m \to + \infty}  r^{-m}\,\sum_{Y\,\in \,  V_m ,\, Y  \underset{m }{\sim}  X} \,\left (u_{\mid V_m}(X)-u_{\mid V_m}(Y)\right ) < + \infty $$

 \noindent We will set:

 $$\partial_n\,u(X) = \displaystyle \lim_{m \to + \infty}  r^{-m}\,\sum_{Y\,\in \,  V_m ,\, Y  \underset{m }{\sim}  X} \,\left (u_{\mid V_m}(X)-u_{\mid V_m}(Y)\right ) < + \infty $$

 \end{definition}

 \vskip 1cm

 \begin{definition}

\noindent For any natural integer~$m$, any~$X$ of~$V_m$, and any continuous function~$u$ on~$\mathfrak {ST} $, we will say that~$u$ admits a normal derivative in~$X$, denoted by~$\partial_n\,u(X)$, if:
 $$ \displaystyle \lim_{k \to + \infty}  r^{-k}\,\sum_{Y\,\in \,  V_k ,\, Y  \underset{k }{\sim}  X} \,\left (u_{\mid V_k} (X)-u_{\mid V_k}(Y)\right ) < + \infty $$

 \noindent We will set:

 $$\partial_n\,u(X) = \displaystyle \lim_{k \to + \infty}  r^{-k}\,\sum_{Y\,\in \,  V_k ,\, Y  \underset{k }{\sim}  X} \,\left (u_{\mid V_k}(X)-u_{\mid V_k}(Y)\right ) < + \infty $$

 \end{definition}

 \vskip 1cm
 \begin{remark} One can thus extend the definition of the normal derivative of~$u$ to~$\mathfrak {ST} $.
 \end{remark}

 \vskip 1cm

 \begin{theorem}

\noindent Let~$u$ be in~\mbox{$\text{dom}\,\Delta$}. The, for any~$X$ of~$\mathfrak {ST} $,~$\partial_n\,u(X)$ exists. Moreover, for any~$v$ of~\mbox{$\text{dom}\,\cal E$}, et any natural integer~$m$, the Gauss-Green formula writes:

$$  {\cal E}(u,v) = -\displaystyle \int_{{\cal D} \left ({\mathfrak {ST}}\right ) } \left (\Delta u \right) \,v \,d\mu + \displaystyle \sum_{V_0} v\,\partial_n\,u  $$
 \end{theorem}
 \vskip 1cm

\section{Spectrum of the Laplacian}

In the following, let~$u$ be in~$\text{dom}\, \Delta$. We will apply the \emph{\textbf{spectral decimation method}} developed by~R.~S.~Strichartz \cite{StrichartzLivre2006}, in the spirit of the works of M.~Fukushima et T.~Shima \cite{Fukushima1992}. In order to determine the eigenvalues of the Laplacian~$\Delta\, u$ built in the above, we concentrate first on the eigenvalues~$\left (-{\lambda_m}\right)_{m\in\N}$ of the sequence of graph Laplacians~$\left (\Delta_m \,u\right)_{m\in\N}$, built on the discrete sequence of graphs~$\left ({\mathfrak {ST}}_m\right)_{m\in\N}$. For any natural integer~$m$, the restrictions of the eigenfunctions of the continuous Laplacian~$\Delta\,u$ to the graph~${\mathfrak {ST}}_m$ are, also, eigenfunctions of the Laplacian~$\Delta_m$, which leads to recurrence relations between the eigenvalues of order~$m$ and~$m+1$.

\vskip 1cm

We thus aim at determining the solutions of the eigenvalue equation:

$$-\Delta\,u=\lambda\,u \quad  \text{on }{ \mathfrak {ST}}$$

\noindent as limits, when the integer~$m$ tends towards infinity, of the solutions of:

$$-\Delta_m\,u=\lambda_m\,u \quad  \text{on }V_m\setminus V_0$$

\noindent We will call them Dirichlet eigenvalues (resp. Neumann eigenvalues) if:
\[u_{\mid \partial{ \mathfrak {ST}}}=0 \quad \left( \text{resp.} \quad {\partial_n u}_{\mid \partial{ \mathfrak {ST}}}=0\right) \]

\vskip 1cm
\noindent Given a strictly positive integer~$m$, let us consider a~$(m-1)-$polyhedron cell, with boundary vertices~$X_0$,~$X_1$,~$X_2,X_3$.\\
\noindent We denote by~$Y_1$,~$Y_2$,~$Y_3$,~$Y_4$,~$Y_5$,~$Y_6$ the points of $V_m\setminus V_{m-1}$ such that:

 \begin{enumerate}

 \item[\emph{i}.]~$Y_1$ is between~$X_0$ and~$X_1$ ;

 \item[\emph{ii}.]~$Y_2$ is between~$X_1$ and~$X_2$ ;

 \item[\emph{iii}.]~$Y_3$ is between~$X_0$ and~$X_2$ ;

 \item[\emph{iv}.]~$Y_4$ is between~$X_0$ and~$X_3$ ;

 \item[\emph{v}.]~$Y_5$ is between~$X_1$ and~$X_3$ ;

 \item[\emph{vi}.]~$Y_6$ is between~$X_2$ and~$X_3$.

 \end{enumerate}
 \vskip 1cm

\noindent  The discrete equation on $\mathfrak{ST}$ leads to the following system:

\begin{align*}
\left(6-\lambda_m \right)u(Y_1)=u(X_0)+u(X_1)+u(Y_2)+u(Y_3)+u(Y_4)+u(Y_5)\\
\left(6-\lambda_m \right)u(Y_2)=u(X_1)+u(X_2)+u(Y_1)+u(Y_3)+u(Y_5)+u(Y_6)\\
\left(6-\lambda_m \right)u(Y_3)=u(X_0)+u(X_2)+u(Y_1)+u(Y_2)+u(Y_4)+u(Y_6)\\
\left(6-\lambda_m \right)u(Y_4)=u(X_0)+u(X_3)+u(Y_1)+u(Y_3)+u(Y_5)+u(Y_6)\\
\left(6-\lambda_m \right)u(Y_5)=u(X_1)+u(X_3)+u(Y_1)+u(Y_2)+u(Y_4)+u(Y_6)\\
\left(6-\lambda_m \right)u(Y_6)=u(X_2)+u(X_3)+u(Y_2)+u(Y_3)+u(Y_4)+u(Y_5)\\
\end{align*}

\noindent Under matricial form, this writes:
\[ \mathbf{x}=\mathbf{A}_m^{-1}\mathbf{b}\]

\noindent with
\[
\mathbf{A}_m =
\begin{pmatrix}
6-\lambda_m & -1 & -1 & -1 & -1 & 0 \\
-1 & 6-\lambda_m & -1 & 0 & -1 & -1 \\
-1 & -1 & 6-\lambda_m & -1 & 0 & -1 \\
-1 & 0 & -1 & 6-\lambda_m & -1 & -1 \\
-1 & -1 & 0 & -1 & 6-\lambda_m & -1 \\
0 & -1 & -1 & -1 & -1 & 6-\lambda_m \\
\end{pmatrix}\]

\[
\mathbf{x} =
\begin{pmatrix}
u(Y_1) \\
u(Y_2) \\
u(Y_3) \\
u(Y_4) \\
u(Y_5) \\
u(Y_6) \\
\end{pmatrix}\]

\[
\mathbf{b} =
\begin{pmatrix}
u(X_0)+u(X_1)\\
u(X_1)+u(X_2)\\
u(X_0)+u(X_2)\\
u(X_0)+u(X_3)\\
u(X_1)+u(X_3)\\
u(X_2)+u(X_3)\\
\end{pmatrix}\]

\noindent By assuming~$\lambda_m \neq \{2,6\}$, one gets:

\begin{align*}
u(Y_1)=\frac{\left(4-\lambda_m\right)(u(X_0)+u(X_1))+2(u(X_2)+u(X_3))}{\left(2-\lambda_m\right) \left(6-\lambda_m\right)}\\
u(Y_2)=\frac{\left(4-\lambda_m\right)(u(X_1)+u(X_2))+2(u(X_0)+u(X_3))}{\left(2-\lambda_m\right) \left(6-\lambda_m\right)}\\
u(Y_3)=\frac{\left(4-\lambda_m\right)(u(X_0)+u(X_2))+2(u(X_1)+u(X_3))}{\left(2-\lambda_m\right) \left(6-\lambda_m\right)})\\
u(Y_4)=\frac{\left(4-\lambda_m\right)(u(X_0)+u(X_3))+2(u(X_1)+u(X_2))}{\left(2-\lambda_m\right) \left(6-\lambda_m\right)}\\
u(Y_5)=\frac{\left(4-\lambda_m\right)(u(X_1)+u(X_3))+2(u(X_0)+u(X_2))}{\left(2-\lambda_m\right) \left(6-\lambda_m\right)}\\
u(Y_6)=\frac{\left(4-\lambda_m\right)(u(X_2)+u(X_3))+2(u(X_0)+u(X_1))}{\left(2-\lambda_m\right) \left(6-\lambda_m\right)}\\
\end{align*}

\noindent Let us now compare the $\lambda_{m-1}-$eigenvalues on~$V_{m-1}$, and the $\lambda_{m}-$eigenvalues on~$V_{m}$. To this purpose, we fix~\mbox{$X_0 \,\in V_{m-1} \setminus V_0$}.\\
 \noindent One has to bear in mind that~$X_0$ also belongs to a~$(m-1)$-cell, with boundary points $X_0$,~$X'_1$,~$X'_2$,~$X'_3$ and interior points~$Z_1$,~$Z_2$,~$Z_3$,~$Z_4$,
 ~ $Z_5$,~$Z_6$.\\
\noindent Thus:

\[ \left(6-\lambda_{m-1} \right)u(X_0)=u(X'_1)+u(X'_2)+u(X'_3)+u(X_1)+u(X_2)+u(X_3) \]

\noindent and:

\[ \left(6-\lambda_{m} \right)u(X_0)=u(Y_1)+u(Y_3)+u(Y_4)+u(Z_1)+u(Z_3)+u(Z_4) \]

\begin{align*}
u(Y_1)+u(Y_3)+u(Y_4)=\frac{\left(8-\lambda_m\right)(u(X_1)+u(X_2)+u(X_3))+3\left(4-\lambda_m\right)u(X_0)}{\left(2-\lambda_m\right) \left(6-\lambda_m\right)}\\
u(Z_1)+u(Z_3)+u(Z_4)=\frac{\left(8-\lambda_m\right)(u(X'_1)+u(X'_2)+u(X'_3))+3\left(4-\lambda_m\right)u(X_0)}{\left(2-\lambda_m\right) \left(6-\lambda_m\right)}
\end{align*}

\noindent By adding member to member, one obtains:

\[
\left(6-\lambda_{m}\right)u(X_0)=\frac{\left(8-\lambda_m\right)\left(6-\lambda_{m-1} \right)+6\left(4-\lambda_m\right)}{\left(2-\lambda_m\right) \left(6-\lambda_m\right)} u(X_0)\\
\]

\noindent We record one more forbidden eigenvalue $\lambda_m=8$, else $\lambda_m$ is independent of $\lambda_{m-1}$. One has:

\[
\left(6-\lambda_{m} \right)^2\left(2-\lambda_m\right)
-6\left(4-\lambda_m\right)=\left(8-\lambda_m\right)\left(6-\lambda_{m-1} \right)\\
\]

\noindent Finally:

\[
\lambda_{m-1}=\lambda_{m}(6-\lambda_{m})
\]

\noindent One may solve:

\[
\lambda_{m}=3\pm \sqrt{9-\lambda_{m-1}}
\]

\noindent Let us introduce:

\[\lambda = 2 \, \displaystyle \lim_{m\rightarrow\infty}6^m \lambda_m\]

\noindent One may note that the limit exists, since, when~$x$ is close to~0:

$$ 3 - \sqrt{9-x}=\displaystyle \frac{x}{6}+{\cal O}(x^2)$$

\noindent Let us now look the Dirichlet eigenvalues and eigenfunctions
\vskip 1cm

\begin{enumerate}

\item \underline{First case :~$m=1$.} \\

The Tetrahedron with its ten vertices can be seen in the following figures:

 \begin{figure}[h!]
 \center{\psfig{height=8cm,width=10cm,angle=0,file=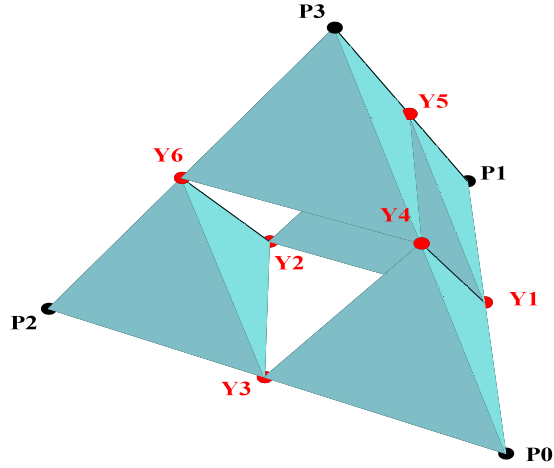}}\\
\caption{The tetrahedron after the first iteration.}
\end{figure}

\noindent Let us look for the kernel of the matrix $\mathbf{A}_1$ in the case where~$\lambda_1\, \in \,\{2,6,8\}$.\\
\noindent For~$\lambda_1=2$, we find the one dimensional Dirichlet eigenspace
$$V^1_2=\text{Vect}\, \left \lbrace (1,1,1,1,1,1)\right \rbrace $$
\noindent For~$\lambda_1=8$, we find the two dimensional Dirichlet eigenspace

$$V^1_8=\text{Vect}\, \left \lbrace (1,-1,0,-1,0,1),(0,-1,1,-1,1,0)\right \rbrace$$

\noindent For~$\lambda_1=6$, we find the three dimensional Dirichlet eigenspace

$$V^1_6=\text{Vect}\, \left \lbrace (-1,0,0,0,0,1),(0,0,-1,0,1,0),(0,-1,0,1,0,0)\right \rbrace$$

\noindent One can easily check that:
$$\#  V_0 =6$$

\noindent Thus, the spectrum is complete.\

\vskip 1cm

\item \underline{Second case :~$m=2$.} \\

\noindent Let us now move to the~$m=2$ case.

 \begin{figure}[h!]
 \center{\psfig{height=8cm,width=10cm,angle=0,file=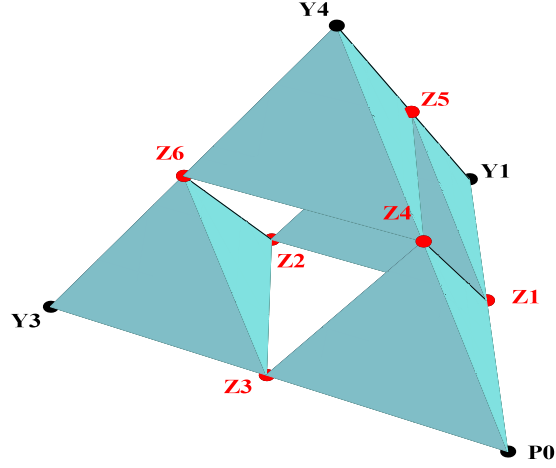}}\\

\caption{The cell $F_1 (V_1)$}
\end{figure}

\noindent Let us denote by $Z^i_j:=f_i(Y_j)$ the points of $V_2\setminus V_1$ that belongs to the cell $f_i(V_0)$.\\

\noindent One has to solve the following system, taking into account the Dirichlet boundary conditions\\~\mbox{$\left (u(X_0)=u(X_1)=u(X_2)=u(X_3)=0\right)$}:

\begin{align*}
\left(6-\lambda_m \right)u(Z^1_1)=u(X_0)+u(Y_1)+u(Z^1_2)+u(Z^1_3)+u(Z^1_4)+u(Z^1_5)\\
\left(6-\lambda_m \right)u(Z^1_2)=u(Y_1)+u(Y_3)+u(Z^1_1)+u(Z^1_3)+u(Z^1_5)+u(Z^1_6)\\
\left(6-\lambda_m \right)u(Z^1_3)=u(X_0)+u(Y_3)+u(Z^1_1)+u(Z^1_2)+u(Z^1_4)+u(Z^1_6)\\
\left(6-\lambda_m \right)u(Z^1_4)=u(X_0)+u(Y_4)+u(Z^1_1)+u(Z^1_3)+u(Z^1_5)+u(Z^1_6)\\
\left(6-\lambda_m \right)u(Z^1_5)=u(Y_1)+u(Y_4)+u(Z^1_1)+u(Z^1_2)+u(Z^1_4)+u(Z^1_6)\\
\left(6-\lambda_m \right)u(Z^1_6)=u(Y_3)+u(Y_4)+u(Z^1_2)+u(Z^1_3)+u(Z^1_4)+u(Z^1_5)\\
\end{align*}

\begin{align*}
\left(6-\lambda_m
\right)u(Z^2_1)=u(X_1)+u(Y_1)+u(Z^2_2)+u(Z^2_3)+u(Z^2_4)+u(Z^2_5)\\
\left(6-\lambda_m \right)u(Z^2_2)=u(X_1)+u(Y_2)+u(Z^2_1)+u(Z^2_3)+u(Z^2_5)+u(Z^2_6)\\
\left(6-\lambda_m \right)u(Z^2_3)=u(Y_1)+u(Y_2)+u(Z^2_1)+u(Z^2_2)+u(Z^2_4)+u(Z^2_6)\\
\left(6-\lambda_m \right)u(Z^2_4)=u(Y_1)+u(Y_5)+u(Z^2_1)+u(Z^2_3)+u(Z^2_5)+u(Z^2_6)\\
\left(6-\lambda_m \right)u(Z^2_5)=u(X_1)+u(Y_5)+u(Z^2_1)+u(Z^2_2)+u(Z^2_4)+u(Z^2_6)\\
\left(6-\lambda_m \right)u(Z^2_6)=u(Y_2)+u(Y_5)+u(Z^2_2)+u(Z^2_3)+u(Z^2_4)+u(Z^2_5)\\
\end{align*}

\begin{align*}
\left(6-\lambda_m \right)u(Z^3_1)=u(Y_2)+u(Y_3)+u(Z^3_2)+u(Z^3_3)+u(Z^3_4)+u(Z^3_5)\\
\left(6-\lambda_m \right)u(Z^3_2)=u(X_2)+u(Y_2)+u(Z^3_1)+u(Z^3_3)+u(Z^3_5)+u(Z^3_6)\\
\left(6-\lambda_m \right)u(Z^3_3)=u(X_2)+u(Y_3)+u(Z^3_1)+u(Z^3_2)+u(Z^3_4)+u(Z^3_6)\\
\left(6-\lambda_m \right)u(Z^3_4)=u(Y_3)+u(Y_6)+u(Z^3_1)+u(Z^3_3)+u(Z^3_5)+u(Z^3_6)\\
\left(6-\lambda_m \right)u(Z^3_5)=u(Y_2)+u(Y_6)+u(Z^3_1)+u(Z^3_2)+u(Z^3_4)+u(Z^3_6)\\
\left(6-\lambda_m \right)u(Z^3_6)=u(X_2)+u(Y_6)+u(Z^3_2)+u(Z^3_3)+u(Z^3_4)+u(Z^3_5)\\
\end{align*}

\begin{align*}
\left(6-\lambda_m \right)u(Z^4_1)=u(Y_4)+u(Y_5)+u(Z^4_2)+u(Z^4_3)+u(Z^4_4)+u(Z^4_5)\\
\left(6-\lambda_m \right)u(Z^4_2)=u(Y_5)+u(Y_6)+u(Z^4_1)+u(Z^4_3)+u(Z^4_5)+u(Z^4_6)\\
\left(6-\lambda_m \right)u(Z^4_3)=u(Y_4)+u(Y_6)+u(Z^4_1)+u(Z^4_2)+u(Z^4_4)+u(Z^4_6)\\
\left(6-\lambda_m \right)u(Z^4_4)=u(X_3)+u(Y_4)+u(Z^4_1)+u(Z^4_3)+u(Z^4_5)+u(Z^4_6)\\
\left(6-\lambda_m \right)u(Z^4_5)=u(X_3)+u(Y_5)+u(Z^4_1)+u(Z^4_2)+u(Z^4_4)+u(Z^4_6)\\
\left(6-\lambda_m \right)u(Z^4_6)=u(X_3)+u(Y_6)+u(Z^4_2)+u(Z^4_3)+u(Z^4_4)+u(Z^4_5)\\
\end{align*}

\begin{align*}
\left(6-\lambda_m \right)u(Y_1)=u(Z^1_1)+u(Z^1_2)+u(Z^1_5)+u(Z^2_1)+u(Z^2_3)+u(Z^2_4)\\
\left(6-\lambda_m \right)u(Y_2)=u(Z^2_2)+u(Z^2_3)+u(Z^2_6)+u(Z^3_1)+u(Z^3_2)+u(Z^3_5)\\
\left(6-\lambda_m \right)u(Y_3)=u(Z^1_2)+u(Z^1_3)+u(Z^1_6)+u(Z^3_1)+u(Z^3_3)+u(Z^3_4)\\
\left(6-\lambda_m \right)u(Y_4)=u(Z^1_4)+u(Z^1_5)+u(Z^1_6)+u(Z^4_1)+u(Z^4_3)+u(Z^4_4)\\
\left(6-\lambda_m \right)u(Y_5)=u(Z^2_4)+u(Z^2_5)+u(Z^2_6)+u(Z^4_1)+u(Z^4_2)+u(Z^4_5)\\
\left(6-\lambda_m \right)u(Y_6)=u(Z^3_4)+u(Z^3_5)+u(Z^3_6)+u(Z^4_2)+u(Z^4_3)+u(Z^4_6)\\
\end{align*}

\noindent The system can be written as $\mathbf{A}_2\mathbf{x}=\mathbf{0}$ and we look for the kernel of $\mathbf{A}_2$ for $\lambda_2 \in \{2,6,8\}$.
\newline We found that there is no eigenfunction for $\lambda_2=2$.

\noindent For $\lambda_2=6$, the eigenspace is a six dimensional one, the basis vectors of which are:

\[
\left(
\begin{array}{c|c|c|c|c|c}
 1 & 0 & -1 & 1 & 1 & 0 \\
 0 & 0 & 0 & -1 & 0 & -1 \\
 0 & 1 & 0 & 0 & -1 & 1 \\
 0 & 0 & 0 & 1 & 0 & 1 \\
 0 & -1 & 0 & 0 & 1 & -1 \\
 -1 & 0 & 1 & -1 & -1 & 0 \\
 0 & 0 & 1 & 0 & -1 & 0 \\
 1 & -1 & 0 & 0 & 1 & -1 \\
 0 & 0 & 0 & 0 & 0 & 1 \\
 -1 & 1 & 0 & 0 & -1 & 1 \\
 0 & 0 & 0 & 0 & 0 & -1 \\
 0 & 0 & -1 & 0 & 1 & 0 \\
 0 & 0 & 0 & 1 & 0 & 0 \\
 -1 & 1 & 1 & -1 & -1 & 0 \\
 0 & 0 & 0 & 0 & 1 & 0 \\
 1 & -1 & -1 & 1 & 1 & 0 \\
 0 & 0 & 0 & 0 & -1 & 0 \\
 0 & 0 & 0 & -1 & 0 & 0 \\
 1 & 0 & 0 & 0 & 0 & 0 \\
 0 & 0 & 1 & 0 & 0 & 0 \\
 0 & 1 & 0 & 0 & 0 & 0 \\
 0 & 0 & -1 & 0 & 0 & 0 \\
 0 & -1 & 0 & 0 & 0 & 0 \\
 -1 & 0 & 0 & 0 & 0 & 0 \\
 0 & 0 & 0 & 0 & 0 & 0 \\
 0 & 0 & 0 & 0 & 0 & 0 \\
 0 & 0 & 0 & 0 & 0 & 0 \\
 0 & 0 & 0 & 0 & 0 & 0 \\
 0 & 0 & 0 & 0 & 0 & 0 \\
 0 & 0 & 0 & 0 & 0 & 0 \\
\end{array}
\right)
\]

\noindent For~$\lambda_2=8$, the eigenspace has dimension 14, and is generated by:
\[
\left(
\begin{array}{c|c|c|c|c|c|c|c|c|c|c|c|c|c}
 0 & 0 & 1 & 1 & 0 & -1 & 0 & 0 & 0 & 0 & 0 & 0 & 1 & 0 \\
 0 & 0 & -1 & -1 & 0 & -1 & 0 & 0 & 0 & 0 & 0 & 0 & -1 & -1 \\
 0 & 0 & 1 & -1 & 0 & 1 & 0 & 0 & 0 & 0 & 0 & 0 & 0 & 1 \\
 0 & 0 & -2 & 0 & 0 & 0 & 0 & 0 & 0 & 0 & 0 & 0 & -1 & -1 \\
 0 & 0 & 0 & 0 & 0 & 0 & 0 & 0 & 0 & 0 & 0 & 0 & 0 & 1 \\
 0 & 0 & 0 & 0 & 0 & 0 & 0 & 0 & 0 & 0 & 0 & 0 & 1 & 0 \\
 0 & 1 & 0 & 0 & 1 & -1 & 0 & 0 & 0 & 0 & 1 & 0 & 0 & 0 \\
 0 & -1 & 0 & 0 & -1 & 1 & 0 & 0 & 0 & 0 & -1 & -1 & 0 & 0 \\
 0 & 1 & 0 & 0 & -1 & -1 & 0 & 0 & 0 & 0 & 0 & 1 & 0 & 0 \\
 0 & -2 & 0 & 0 & 0 & 0 & 0 & 0 & 0 & 0 & -1 & -1 & 0 & 0 \\
 0 & 0 & 0 & 0 & 0 & 0 & 0 & 0 & 0 & 0 & 0 & 1 & 0 & 0 \\
 0 & 0 & 0 & 0 & 0 & 0 & 0 & 0 & 0 & 0 & 1 & 0 & 0 & 0 \\
 1 & 0 & 0 & -1 & -1 & 0 & 0 & 0 & 1 & 0 & 0 & 0 & 0 & 0 \\
 -1 & 0 & 0 & 1 & -1 & 0 & 0 & 0 & -1 & -1 & 0 & 0 & 0 & 0 \\
 1 & 0 & 0 & -1 & 1 & 0 & 0 & 0 & 0 & 1 & 0 & 0 & 0 & 0 \\
 -2 & 0 & 0 & 0 & 0 & 0 & 0 & 0 & -1 & -1 & 0 & 0 & 0 & 0 \\
 0 & 0 & 0 & 0 & 0 & 0 & 0 & 0 & 0 & 1 & 0 & 0 & 0 & 0 \\
 0 & 0 & 0 & 0 & 0 & 0 & 0 & 0 & 1 & 0 & 0 & 0 & 0 & 0 \\
 1 & -1 & -1 & 0 & 0 & 0 & 1 & 0 & 0 & 0 & 0 & 0 & 0 & 0 \\
 -1 & -1 & 1 & 0 & 0 & 0 & -1 & -1 & 0 & 0 & 0 & 0 & 0 & 0 \\
 -1 & 1 & -1 & 0 & 0 & 0 & 0 & 1 & 0 & 0 & 0 & 0 & 0 & 0 \\
 0 & 0 & 0 & 0 & 0 & 0 & -1 & -1 & 0 & 0 & 0 & 0 & 0 & 0 \\
 0 & 0 & 0 & 0 & 0 & 0 & 0 & 1 & 0 & 0 & 0 & 0 & 0 & 0 \\
 0 & 0 & 0 & 0 & 0 & 0 & 1 & 0 & 0 & 0 & 0 & 0 & 0 & 0 \\
 0 & 0 & 0 & 0 & 0 & 2 & 0 & 0 & 0 & 0 & 0 & 0 & 0 & 0 \\
 0 & 0 & 0 & 0 & 2 & 0 & 0 & 0 & 0 & 0 & 0 & 0 & 0 & 0 \\
 0 & 0 & 0 & 2 & 0 & 0 & 0 & 0 & 0 & 0 & 0 & 0 & 0 & 0 \\
 0 & 0 & 2 & 0 & 0 & 0 & 0 & 0 & 0 & 0 & 0 & 0 & 0 & 0 \\
 0 & 2 & 0 & 0 & 0 & 0 & 0 & 0 & 0 & 0 & 0 & 0 & 0 & 0 \\
 2 & 0 & 0 & 0 & 0 & 0 & 0 & 0 & 0 & 0 & 0 & 0 & 0 & 0 \\
\end{array}
\right)
\]

\noindent From~$\lambda_1=2$, the spectral decimation leads to:

$$\lambda_{2}=3- \sqrt{7} \quad \text{and} \quad \lambda_{2}=3+ \sqrt{7}$$

\noindent Each of these eigenvalues has multiplicity~$1$.\\

\noindent From~$\lambda_1=6$, the spectral decimation leads to:
$$\lambda_{2}=3- \sqrt{3} \quad \text{and} \quad \lambda_{2}=3+ \sqrt{3}$$

\noindent Each of these eigenvalues has multiplicity~$3$.\\

\noindent From~$\lambda_1=8$, the spectral decimation leads to:
$$\lambda_{2}=4$$
\noindent with multiplicity~$2$ (one may note that~$2$ is not a Dirichlet eigenvalue for $m=2$).\\

\noindent One can easily check that:

$$\# V_2\setminus V_0 =30=6+14+2\times 1+3\times 2 + 1\times 2=30$$

\noindent Thus,the spectrum is complete.

\end{enumerate}

\vskip 1cm

\noindent Let us now go back to the general case. Given a strictly positive integer~$m$, let us introduce the respective multiplicities~$M_m(6)$ and $M_m(8)$ of the eigenvalues $\lambda_m=6$ and $\lambda_m=8$.\\

\noindent One can easily check by induction that:

$$\# V_m\setminus V_0 =2(4^m -1 )$$

\noindent and:
$$M_m(8)=4^{m} - 2$$
\noindent (we here refer to~\cite{Shima}).\\
\noindent One also has:

$$\# V_{m-1}\setminus V_0 =2(4^{m-1} -1 )$$
\noindent and:
$$M_{m-1}(8)=4^{m-1} - 2$$

\noindent There are thus~\mbox{ $2(4^{m-1} -1)-(4^{m-1} - 2)=4^{m-1}$} continued eigenvalues (the ones obtained by means of the spectral decimation), which correspond to a space of eigenfunctions, the dimension of which is:

$$ (4^{m-1} - 2) + 2\times 4^{m-1}$$

\noindent This leads to:
$$M_m(6)=2(4^m -1 )-(4^{m} - 2) -\left ((4^{m-1} - 2) + 2\times 4^{m-1}\right)=4^m$$

\section{Metric - Towards spectral asymptotics}

\begin{definition}\textbf{Effective resistance metric, on~${\mathfrak {ST}}$}\\
\noindent Given two points~$(X,Y)$ of~\mbox{${\mathfrak {ST}}^2$}, let us introduce the \textbf{effective resistance metric between~$X$ and~$Y$}:

$$ R_{\mathfrak {ST}}(X,Y)= \left \lbrace \min_{\left \lbrace u \, |\, u(X)=0 , u(Y)=1\right \rbrace } \,{\cal E}(u )  \right \rbrace^{-1}$$

\noindent In an equivalent way,~\mbox{$ R_{\mathfrak {ST}}(X,Y)$} can be defined as the minimum value of the real numbers~$R$ such that, for any function~$u$ of~\mbox{$\text{dom} \,\Delta$}:

$$ \left |u (X)-u (Y)\right |^2 \leq R\,  {\cal E}(u)  $$

\end{definition}

\vskip 1cm

\begin{definition}\textbf{Metric, on the Sierpi\'{n}ski Tetrahedron~${\mathfrak {ST}}$}\\
\noindent Let us define, on the Sierpi\'{n}ski Tetrahedron~${\mathfrak {ST}}$, the distance~\mbox{$d_{ \mathfrak {ST} }$} such that, for any pair of points~$(X,Y)$
 of~\mbox{$ {\mathfrak {ST}}^2 $}:

$$ d_{{\mathfrak {ST}}}(X,Y)= \left \lbrace \min_{\left \lbrace u \, |\, u(X)=0 , u(Y)=1\right \rbrace } \,{\cal E}(u,u)  \right \rbrace^{-1}$$
\end{definition}

\vskip 1cm

\begin{remark}

One may note that the minimum

$$\min_{\left \lbrace u \, |\, u(X)=0 , u(Y)=1\right \rbrace } \,{\cal E}(u )   $$

\noindent is reached for~$u$ being harmonic on the complement set, on~\mbox{$\mathfrak {ST}$}, of the set
$$ \left \lbrace X \right \rbrace \cup  \left \lbrace Y \right \rbrace  $$

\noindent (One might bear in mind that, due to its definition, a harmonic function~$u$ on~
\mbox{$\mathfrak {ST}$} minimizes the sequence of energies~\mbox{$\left ({\cal E}_{ {\mathfrak {ST}}_m }  (u,u )\right)_{m\in\N}$}.\\

\end{remark}

\vskip 1cm

\begin{definition}\textbf{Dimension of the Sierpi\'{n}ski Tetrahedron~${\mathfrak {ST}}$, in the resistance metric}\\

\noindent The \textbf{dimension of the Sierpi\'{n}ski Tetrahedron~${\mathfrak {ST}}$}, in the resistance metrics, is the strictly positive number~\mbox{$d_{\mathfrak {ST} }$} such that, given a strictly positive real number~$r$, and a point~$X\,\in\,{\mathfrak {ST}} $, for the~\mbox{$X-$centered} ball of radius~$r$, denoted by~\mbox{${\cal B}_r(X)$}:

$$ \mu \left ({\cal B}_r(X) \right) =r^{d_{\mathfrak {ST}}}$$

\end{definition}

\vskip 1cm

\begin{pte}

\noindent Given a natural integer~$m$, and two points~$(X,Y)$ of~\mbox{$ {\mathfrak {ST}}^2 $} such that~\mbox{$X  \underset{m }{\sim}  Y $}:

 $$ \displaystyle \min_{\left \lbrace u \, |\, u(X)=0 , u(Y)=1\right \rbrace } \,{\cal E}(u ) \lesssim r^m=\left (\displaystyle \frac{2}{3} \right)^m$$

 \noindent which also corresponds to the order of the diameter of~$m-$polyhedra cells. \\

\noindent Since the Sierpi\'{n}ski tetrahedron~$\mathfrak {ST}$ is obtained from the initial regular~\mbox{$3-$simplex} by means of four contractions, the ratio which is equal to~$\displaystyle \frac{1}{2}$, let us look for a real number~$\beta_{\mathfrak {ST}}$ such that:

 $$\left (\displaystyle \frac{1}{2} \right)^{m\,\beta_{\mathfrak {ST}}}=\left (\displaystyle \frac{2}{3} \right)^m$$

 \noindent One obtains:

 $$ \beta_{\mathfrak {ST}} = \displaystyle \frac{\ln \frac{3}{2}}{\ln 2}  $$

 \noindent Let us denote by~$\mu$ the standard measure on~$\mathfrak {ST}$ which assigns measure~\mbox{$\displaystyle \frac{1}{4^m}$ } to each $m-$polyhedron cell. Let us now look for a real number~$d_{\mathfrak {ST}}$ such that:

$$\left (\displaystyle \frac{2}{3 }\right)^{m \,d_{\mathfrak {ST}}}= \displaystyle \frac{1}{4^m}$$

 \noindent One obtains:

 $$d_{\mathfrak {ST} }=\displaystyle \frac{ \ln   \frac{3}{2}}{\ln 4}$$

 \noindent   Given s strictly positive real number~$r$, and a point~$X\,\in\,{\mathfrak {ST}} $, one has then the following estimate, for the~\mbox{$X-$centered} ball of radius~$r$, denoted by~\mbox{${\cal B}_r(X)$}:

$$ \mu \left ({\cal B}_r(X) \right) =r^{d_{\mathfrak {ST}}}$$

\end{pte}

\vskip 1cm

\begin{definition}\textbf{Eigenvalue counting function}\\

\noindent We introduce the eigenvalue counting function~${\cal N}^{\mathfrak {ST}}$ such that, for any real number~$x$:

  $${\cal N}^{\mathfrak {ST}}(x) = \#\, \left \lbrace \lambda\,\text{eigenvalue of $-\Delta$} \, : \quad \lambda \leq  x \right \rbrace $$

\end{definition}

\vskip 1cm

\begin{pte}
\noindent The existing results of~J.~Kigami~\cite{Kigami1998} and~R.~S.~Strichartz~\cite{Strichartz2012} lead to the modified Weyl formula, when~$x$ tends towards infinity:

  $${\cal N}^{\mathfrak {ST}}(x)= G(x)\, x^{\alpha_{\mathfrak {ST}}}+ {\cal O}(1)$$

\noindent where the exponent~$\alpha_{\mathfrak {ST} }$ is given by:

 $$\alpha_{\mathfrak {ST} }=\displaystyle \frac{ d_{\mathfrak {ST} }}{d_{\mathfrak {ST} }+1}$$

\end{pte}

\end{document}